\documentclass[12pt]{article}
\usepackage{amsmath,amssymb,amscd}
\newcommand{\x}{\times}

\renewcommand{\>}{\rangle}

\renewcommand{\a}{\alpha}
\renewcommand{\b}{\beta}
\renewcommand{\d}{\delta}

\newcommand{\e}{\varepsilon}
\newcommand{\g}{\gamma}
\newcommand{\G}{\Gamma}

\newcommand{\var}{\varphi}
\newcommand{\s}{\sigma}

\renewcommand{\O}{\Omega}

\renewcommand{\i}{\infty}

\newcommand{\mcap}{\mathop{\mathrm{cap}}\nolimits}
\newcommand{\mtop}{\mathop{\mathrm{top}}\nolimits}
\newcommand{\mvol}{\mathop{\mathrm{vol}}\nolimits}
\newcommand{\mop}{\mathop{\mathrm{op}}\nolimits}

\title{The Topological Version of Free Entropy}
\author{Dan Voiculescu\footnote{Supported in part by NSF Grant
DMS-0079945.}\\
Department of Mathematics\\
University of California at Berkeley\\
Berkeley, CA\ \ 94720-3840\\
E-mail: dvv@math.berkeley.edu}
\date{}

\begin{document}
\maketitle

\begin{abstract}
Notions of topological free entropy and of free capacity are
introduced in the $C^*$-algebra context.  Basic properties, basic problems
and connections to potential theory and random matrix theory are
discussed.
\end{abstract}

\setcounter{section}{-1}
\section{Introduction}
\label{sec0}

The game of defining free entropy via microstates \cite{8} can also be
played with norm-microstates:  instead of expectation values one
approximates norms of polynomials.  This leads to an invariant associated
with $n$-tuples of elements of a $C^*$-algebra:  the {\em topological free
entropy}.

A competing notion is the {\em free capacity}, which we define as the
supremum over trace-states of the free entropy.  In the one-variable case
there is a simple relation to the logarithmic capacity of potential theory
\cite{4}.

We prove that the free capacity majorizes the topological free entropy and
this leads to the obvious question about equality.  In the two examples
which we study, the free semicircular $n$-tuple and the universal
$n$-tuple of self-adjoint contraction, the two quantities are equal.

Notions of free entropy dimension in $C^*$-algebras can also be defined
either via a supremum over trace-states or directly from the
norm-microstates.  In the latter case, instead of adapting our initial
definition \cite{8}, which was inspired by the Minkowski content, we found
it more convenient to use the recent simpler approach of Jung \cite{3}.

We began exploring topological free entropy a few years ago, but abandoned
it soon after facing the difficulties of the semicircular $n$-tuple
example.  The impetus to return to this work was provided by the recent
breakthrough of U.~Haagerup and S.~Thorbjornsen \cite{2} on random
matrices, which also gave what was essentially needed to answer the
semicircular question.  There is good reason now, we believe, to expect
that topological free entropy will be useful in the $C^*$-algebra context
and that progress will be made on the problems we discuss in this note.

\section{Free Capacity}
\label{sec1}

Let $A$ be a unital $C^*$-algebra generated by $n+m$ self-adjoint elements
$a_1,\dots,a_n,b_1,\dots,b_m$.  Let further $TS(A)$ denote the
trace-states of $A$.  Given $\tau \in TS(A)$ an element $a \in A$ gives
rise to a noncommutative random variable in the tracial $C^*$-probability
space $(A,\tau)$, which we shall denote by $a(\tau)$.  The free entropy
$\chi(a_1(\tau),\dots,a_n(\tau): b_1(\tau),\dots,b_m(\tau))$ of the
$a_i(\tau)$'s in the presence of the $b_j(\tau)$'s \cite{9} will also be
denoted $\chi(a_1,\dots,a_n: b_1,\dots,b_m; \tau)$ (in case $m = 0$ we
shall write $\chi(a_1,\dots,a_n;\tau)$).

We define the {\em free capacity $\kappa(a_1,\dots,a_n: b_1,\dots,b_m)$
of $a_1,\dots,a_n$ in the presence of $b_1,\dots,b_m$ by}:
\[
\kappa(a_1,\dots,a_n: b_1,\dots,b_m) = \sup_{\tau \in TS(A)}
\chi(a_1,\dots,a_n: b_1,\dots,b_m;\tau).
\]

Because of Corollary $1.8$ in \cite{9} $\kappa(a_1,\dots,a_n:
b_1,\dots,b_m)$ does not depend on the choice of $b_1,\dots,b_m$ as long
as $A$ is generated by $a_1,\dots,a_n,b_1,\dots,b_m$.  Hence it is also
appropriate to use {\em the notation $\kappa(a_1,\dots,a_n: A)$ for
$\kappa(a_1,\dots,a_n: b_1,\dots,b_m)$.  If $A$ is already generated by
$a_1,\dots,a_n$ we shall write simply $\kappa(a_1,\dots,a_n)$ for
$\kappa(a_1,\dots,a_n: A)$ and call it the free capacity of
$a_1,\dots,a_n$}.

In case $n = 1$, up to sign and to the additive constant $\theta = \frac
{3}{4} + \frac {1}{2} \log 2\pi$ the free entropy is the logarithmic
energy of the distribution.  Therefore if $\sigma(a)$ is the spectrum of
$a = a^*$, we have
\[
e^{\kappa(a)-\theta} = \mcap(\s(a))
\]
where $\mcap(\cdot)$ denotes the logarithmic capacity (\cite{4}).

We chose to call $\kappa(a_1,\dots,a_n)$, instead of its exponential, the
free capacity, because from an information theory point of view, it is the
largest quantity of ``free information'' which can be realized by
$a_1,\dots,a_n$.  However, as shown for instance by the free analogue of
the entropy-power inequality, \cite{7}, exponentials of free entropies are
also quantities of interest.

In view of the upper semicontinuity of the free entropy with respect to
convergence in distribution and compactness of $TS(A)$ we infer that:
\[
\kappa(a_1,\dots,a_n: b_1,\dots,b_m) = \chi(a_1,\dots,a_n:
b_1,\dots,b_m;\tau)
\]
{\em for some} $\tau \in TS(A)$ (actually an extremal trace-state in view
of \cite{9}).

Also, among the properties of $\kappa(a_1,\dots,a_n)$ which follow easily
from those of $\chi$ we note two more.  Subadditivity of $\kappa$ holds
\begin{eqnarray*}
\kappa(a_1,\dots,a_n,b_1,\dots,b_m) &\le &\kappa(a_1,\dots,a_n:
b_1,\dots,b_m) + \kappa(b_1,\dots,b_m: a_1,\dots,a_n) \\
&\le &\kappa(a_1,\dots,a_n) + \kappa(b_1,\dots,b_m).
\end{eqnarray*}
For a change of variables given by noncommutative power series satisfying
all the conditions of Proposition $3.5$ in \cite{8}, we get:
\begin{eqnarray*}
\kappa(F_1(a_1,\dots,a_n),\dots,F_n(a_1,\dots,a_n)) &\le
&\kappa(a_1,\dots,a_n) \\
& &+ \sup_{\tau \in TS(A)} \log|{\mathcal
J}|((F_1,\dots,F_n))(a_1(\tau),\dots,a_n(\tau)).
\end{eqnarray*}

Finally, a remark which applies not only to free capacity, but also to the
other constructs in this note.  The appropriate $C^*$-algebras for these
considerations are those subject to finiteness conditions, such as having
sufficiently many trace-states.  Of course, given an arbitrary unital
$C^*$-algebra, we can always factor by the two-sided ideal on which all
trace-states vanish and work in the quotient.

\section{Topological Free Entropy}
\label{sec2}

Let $A,a_1,\dots,a_n,b_1,\dots,b_m$ be like in the preceding section. 
Let further 
\[
P_1,\dots,P_r \in {\mathbb
C}\<X_1,\dots,X_n,Y_1,\dots,Y_m\>
\]
be polynomials in the noncommutative
indeterminates $X_1,\dots,X_n,Y_1,\dots,Y_m$ and let $\e > 0$, $k \in
{\mathbb N}$.  We define the {\em norm-microstates (or topological
microstates)}
\[
\G_{\mtop}(a_1,\dots,a_n,b_1,\dots,b_m; k,\e,P_1,\dots,P_r)
\]
to be $n+m$-tuples of self-adjoint $k \x k$ matrices
$(C_1,\dots,C_n,D_1,\dots,D_m) \in ({\mathfrak M}_k^{sa})^{n+m}$ so that
\[
|\|P_j(C_1,\dots,C_n,D_1,\dots,D_m)\| -
\|P_j(a_1,\dots,a_n,b_1,\dots,b_m)\|| \le \e
\]
\[
1 \le j \le r.
\]
We define the norm-microstates for $a_1,\dots,a_n$ in the presence of
$b_1,\dots,b_m$, denoted
\[
\G_{\mtop}(a_1,\dots,a_n: b_1,\dots,b_m; k,\e,P_1,\dots,P_r)
\]
to be the projection via
\[
(C_1,\dots,C_n,D_1,\dots,D_m) \rightarrow (C_1,\dots,C_n)
\]
of
\[
\G_{\mtop}(a_1,\dots,a_n,b_1,\dots,b_m; k,\e,P_1,\dots,P_r)
\]
into $({\mathfrak M}_k^{sa})^n$.

With this change of microstates, the definition of the {\em topological
free entropy} $\chi_{\mtop}$ proceeds exactly like the definition of
$\chi$.  We take
\[
\limsup_{k \rightarrow \i} \left(k^{-2}\log\mvol \G_{\mtop}(\dots) + \frac
{n}{2} \log k\right)
\]
followed by
\[
\inf_{\e > 0} \inf_{r \in {\mathbb N}} \inf_{P_1,\dots,P_r \in {\mathbb
C}\<X_1,\dots,X_n,Y_1,\dots,Y_m\>}.
\]

We call the resulting $\chi_{\mtop}(a_1,\dots,a_n: b_1,\dots,b_m)$ {\em
the topological free entropy of $a_1,\dots,a_n$ in the presence of
$b_1,\dots,b_m$.  In case $m=0$ we write $\chi_{\mtop}(a_1,\dots,a_n)$ and
call it the topological free entropy of $a_1,\dots,a_n$.}

Several basic properties of $\chi_{\mtop}$ are proved along similar lines
to properties of $\chi$ and we leave the details to the reader.

\subsection*{Subadditivity}

\begin{eqnarray*}
\chi_{\mtop}(a_1,\dots,a_n,b_1,\dots,b_m) &\le
&\chi_{\mtop}(a_1,\dots,a_n: b_1,\dots,b_m) \\
& &+ \chi_{\mtop}(b_1,\dots,b_m: a_1,\dots,a_n) \\
&\le &\chi_{\mtop}(a_1,\dots,a_n) + \chi_{\mtop}(b_1,\dots,b_m).
\end{eqnarray*}

\subsection*{Upper Semicontinuity}

Assume $a_1^{(p)},\dots,a_n^{(p)} \in A_p$ converge in
``norm-distribution'' to $a_1,\dots,a_n \in A$ in the sense, that for
every noncommutative polynomial $P \in {\mathbb C}\<X_1,\dots,X_n\>$ we
have
\[
\lim_{p \rightarrow \i} \|P(a_1^{(p)},\dots,a_n^{(p)})\| =
\|P(a_1,\dots,a_n)\|.
\]
Then we have
\[
\limsup_{p \rightarrow \i} \chi_{\mtop}(a_1^{(p)},\dots,a_n^{(p)}) \le
\chi_{\mtop}(a_1,\dots,a_n).
\]

\subsection*{Change of Variables}

Assume there are noncommutative power series $F_1,\dots,F_n,G_1,\dots,G_n$
satisfying the complete set of conditions in Proposition $3.5$ of
\cite{8}.  Then we have
\begin{eqnarray*}
\chi_{\mtop}(F_1(a_1,\dots,a_n),\dots,F_n(a_1,\dots,a_n)) &\le
&\chi_{\mtop}(a_1,\dots,a_n) \\
& &+ n \log\|{\tilde D}F(a_1,\dots,a_n)\|.
\end{eqnarray*}
Here ${\tilde D}F(a_1,\dots,a_n) \in {\mathfrak M}_n \otimes A\ 
{\bar \otimes}\ A^{\mop}$ where the tensor product is the spatial tensor
product.  To define ${\tilde D}F(a_1,\dots,a_n)$ one first notices that in
a formal sense the differential of $F$ as power-series is a $n \x n$ matrix
of partial free difference quotients and that, under the assumptions on
multiradii of convergence in Proposition $3.5$ of \cite{8}, this
differential, after replacing the indeterminates by $a_1,\dots,a_n$ gives
rise to an element in the projective tensor product ${\mathfrak M}_n
\otimes A \otimes_{\pi} A^{\mop}$.  Clearly this element gives then rise
to an element in ${\mathfrak M}_n \otimes A\ {\bar \otimes}\ A^{\mop}$. 
The proof proceeds along the same lines as the proof of Proposition $3.5$
in \cite{8}.

\section{The Inequality $\chi_{\mtop} \le \kappa$}
\label{sec3}

{\em Under the same assumptions about $A,a_1,\dots,a_n,b_1,\dots,b_m$ as
in the previous sections, we shall prove here that}
\[
\chi_{\mtop}(a_1,\dots,a_n: b_1,\dots,b_m) \le \kappa(a_1,\dots,a_n:
b_1,\dots,b_m).
\]

If $\chi_{\mtop}(a_1,\dots,a_n: b_1,\dots,b_m) = -\i$ there is nothing to
prove.  Assume
\[
\chi_{\mtop}(a_1,\dots,a_n: b_1,\dots,b_m) > \a > -\i.
\]
We will show that $\kappa(a_1,\dots,a_n: b_1,\dots,b_m) > \a$.

Let ${\mathbb Q}[i]\<X_1,\dots,X_n,Y_1,\dots,Y_m\>$ be the noncommutative
polynomials with rational complex coefficients, which being a countable
set can be put in a sequence $P_1,P_2,\dots$.  It will be also convenient
to assume $P_1,\dots,P_{m+n}$ are precisely $X_1,\dots,X_n,Y_1,\dots,Y_m$.

It is easy to see that 
\[
\chi_{\mtop}(a_1,\dots,a_n: b_1,\dots,b_m) > \a
\]
is equivalent to the existence of integers $m+n < k_1 < k_2 < \dots$ so
that for some $\a' > \a$ 
\[
k_r^{-2} \log \mvol\ \G_{\mtop}(a_1,\dots,a_n: b_1,\dots,b_m;
P_1,\dots,P_r,k_r,r^{-1}) + \frac {n}{2} \log k_r > \a'.
\]
Note that if
\[
\eta = (C_1,\dots,C_n,D_1,\dots,D_m) \in
\G_{\mtop}(a_1,\dots,a_n,b_1,\dots,b_m; P_1,\dots,P_r,k_r,r^{-1}) = \G(r)
\]
then $\|C_k\| < M$, $\|D_j\| < M$ for some $M$ independent of $r$ and of
the choice of microstate.  Let then $A(n+m)$ be the universal unital
$C^*$-algebra generated by self-adjoint elements
$x_1,\dots,x_n,y_1,\dots,y_m$ of norm $M$, that is the unital full free
product of $n+m$ copies of $C([-M,M])$.  A microstate $\eta \in \G(r)$ as
above defines a $*$-homomorphism $\var[\eta]: A(n+m) \rightarrow
{\mathfrak M}_{k_r}$ so that $\var[\eta](x_k) = C_k$, $\var[\eta](y_j) =
D_j$ and a trace-state $\tau[\eta] \in TS(A(n+m))$
\[
\tau[\eta] = k_r^{-1}Tr_{k_r} \circ \var[\eta].
\]
Similarly there is a $*$-homomorphism $\var: A(n+m) \rightarrow A$ so that
$\var(x_k) = a_k$, $\var(y_j) = b_j$.

On $\O = TS(A(n+m))$ the weak topology is induced by the metric
\[
d(\tau_1,\tau_2) = \sum_{p \ge 1} \sum_{(i_1,\dots,i_p) \in
(\{1,\dots,n+m\})^p} (2M(n+m))^{-p}
|(\tau_1-\tau_2)(z_{i_1}\dots z_{i_p})|
\]
where $z_1,\dots,z_{n+m}$ denotes $x_1,\dots,x_n,y_1,\dots,y_m$.  $\O$ is a
compact metric space and
\[
K_r = \{\tau[\eta] \in \O \mid \eta \in \G(r)\}
\]
are compact subsets because $\eta \rightarrow \tau[\eta]$ is continuous
and $\G(r)$ is compact.  Let further $K \subset \O$ denote the compact
subset $(TS(A)) \circ \var$.

If $\tau \in \O$ is the weak limit of some sequence $(\tau[\eta_{r_s}])_{s
\in {\mathbb N}}$, where $r_1 < r_2 < \dots$ and $\tau[\eta_{r_s}] \in
K_{r_s}$, then $\tau \in K$.  Indeed we have
\begin{eqnarray*}
|\tau(P_j(x_1,\dots,x_n,y_1,\dots,y_m))| &\le &\limsup_{s \rightarrow \i}
\|\var[\eta_{r_s}](P_j(x_1,\dots,x_n,y_1,\dots,y_m))\| \\
&\le &\limsup_{s \rightarrow \i} (r_s^{-1} +
\|P_j(a_1,\dots,a_n,b_1,\dots,b_m)\|) \\
&\le &\|P_j(a_1,\dots,a_n,b_1,\dots,b_m)\| \\
&= &\|\var(P_j(x_1,\dots,x_n,y_1,\dots,y_m))\|
\end{eqnarray*}
which in view of the density of the $P_j$'s in $A(n+m)$ gives $\tau \in
K$.  Given $\e > 0$, in view of compactness of $\O$ there is $L(\e) > 0$
so that for each $r$, $K_r = K_r^{(1)} \cup \dots \cup K_r^{(L(\e))}$
where each compact set $K_r^{(j)}$ has diameter $\le \frac {\e}{2}$.  Let
\[
\G(r,j) = \{\eta \in \G(r) \mid \tau[\eta] \in K_r^{(j)}\}.
\]
We have $\G(r) = \G(r,1) \cup \dots \cup \G(r,L(\e))$.  Let further
\[
\pi_r: \G(r) \rightarrow \G_{\mtop}(a_1,\dots,a_n: b_1,\dots,b_m;
P_1,\dots,P_r,k_r,r^{-1})
\]
denote the projection mapping $(C_1,\dots,C_n,D_1,\dots,D_m)$ to
$(C_1,\dots,C_n)$.  For each $r$ let $\G'(r)$ denote some $\pi_r(\G(r,j))$
such that
\[
\mvol(\pi_r(\G(r,j))) \ge (L(\e))^{-1} \mvol(\G_{\mtop}(a_1,\dots,a_n:
b_1,\dots,b_m; P_1,\dots,P_r,k_r,r^{-1})).
\]
Since $k_r \rightarrow \i$ we will clearly have 
\[
\liminf_{r \rightarrow \i} \left(k_r^{-2} \log \mvol(\G'(r)) + \frac
{n}{2} \log k_r\right) \ge \a'.
\]

Giving $\e$ successively the values $1,1/2,1/3,\dots$ and using the fact
that accumulation points of the sequence of sets $K_r$ are in $K$ we find
there is $\tau \in K$ and there are $r_1 < r_2 < \dots$ such that the
chosen set $K_{r_s}^{(j_s)} \subset K_{r_s}$ for $\e = (2s)^{-1}$ is
$\subset B(\tau,s^{-1})$ the ball of radius $s^{-1}$ and center $\tau$ and
the corresponding set $\G'(r_s)$ is such that 
\[
\liminf_{s \rightarrow \i} \left( k_{r_s}^{-2} \log \mvol(\G'(r_s)) +
\frac {n}{2} \log k_{r_s}\right) \ge \a'.
\]
Since $K_{r_s}^{(j_s)} \subset B(\tau,s^{-1})$ we will have
\[
\G'(r_s) \subset \G_M(a_1,\dots,a_n: b_1,\dots,b_m;k_{r_s},m_s,\e_s;\tau)
\]
for some $m_s\rightarrow \i$ and $\e_s \rightarrow 0$.  Clearly this
implies
\[
\chi(a_1,\dots,a_n: b_1,\dots,b_m;\tau) \ge \a'
\]
and hence
\[
\kappa(a_1,\dots,a_n: b_1,\dots,b_m) > \a.
\]

\section{Semicircular Systems}
\label{sec4}

Let $A$ be the $C^*$-algebra generated by a semicircular system
$S_1,\dots,S_n$ consisting of $n$ freely independent $(0,1)$ semicircular
random variables in some $W^*$-probability space $(M,\tau)$ with a
faithful trace-state.  By a recent result of U.~Haagerup and
S.~Thorbjornsen \cite{2} for any $\e > 0$ and noncommutative polynomials
in $n$ indeterminates $P_1,\dots,P_r$ we have
\[
\lim_{k \rightarrow \i} \g_k(\G_{\mtop}(S_1,\dots,S_n;P_1,\dots,P_r,k,\e))
= 1
\]
where $\g_k$ is the Gaussian probability measure on $({\mathfrak
M}_k^{sa})^n$ with density 
\[
c_k \exp\left( -\frac {k}{2} Tr(A_1^2 + \dots + A_n^2)\right)
\]
where $c_k$ is a normalizing constant.  By Proposition $3.1$ of \cite{1},
$A$ has the Dixmier property hence $[A,A] + {\mathbb C}1$ is norm-dense in
$A$.  In particular there are polynomials
\[
Q_{1,j},\dots,Q_{s,j},R_{1,j},\dots,R_{s,j}
\]
such that for a given $\d >
0$ we have
\[
\left\|S_j^2 - 1 - \sum_{1 \le i \le s}
[Q_{ij}(S_1,\dots,S_n),R_{ij}(S_1,\dots,S_n)]\right\| < \d/2.
\]

Thus for an appropriate choice of polynomials $P_1,\dots,P_r$ and $\e > 0$
we will have for all $k$
\begin{eqnarray*}
\G(k) &= &\G_{\mtop}(S_1,\dots,S_n;P_1,\dots,P_r,k,\e) \\
&\subset &\{(A_j)_{1
\le j \le n} \in ({\mathfrak M}_k^{sa})^n \mid k^{-1}Tr_kA_j^2 < 1 + \d, 1
\le j \le n\}.
\end{eqnarray*}
Hence
\[
\g_k(\G(k)) \le c_k \exp\left( -\frac {k^2}{2}
(1+\d)^2n\right)\mvol(\G(k)).
\]
This gives
\[
k^{-2} \log \mvol(\G(k)) + \frac {n}{2} \log k \ge k^{-2}\log \g_k(\G(k))
- k^{-2} \log c_k + \frac {n}{2} \log k - \frac {n}{2} (1+\d)^2.
\]
Since $\g_k(\G(k)) \rightarrow 1$ as $k \rightarrow \i$ and all this holds
also after decreasing $\e > 0$ and enlarging the set of polynomials, we
infer that 
\[
\chi_{\mtop}(S_1,\dots,S_n) \ge \limsup_{k \rightarrow \i} \left( \frac
{n}{2} \log k - k^{-2} \log c_k - \frac {n}{2} (1+\d)^2\right).
\]
Since $c_k = (2\pi)^{-nk^2/2}(k)^{nk^2/2}$ this gives
\[
\chi_{\mtop}(S_1,\dots,S_n) \ge \frac {n}{2} \log 2\pi - \frac
{n}{2}(1+\d)^2.
\]

Since $\d > 0$ is arbitrary we have that
\[
\chi_{\mtop}(S_1,\dots,S_n) \ge \frac {n}{2} (\log 2\pi - 1).
\]
The right-hand side is precisely $\chi(S_1,\dots,S_n;\tau)$.  So we have
\[
\chi_{\mtop}(S_1,\dots,S_n) \ge \chi(S_1,\dots,S_n;\tau).
\]
On the other hand since $A$ has a unique trace-state 
\[
\chi(S_1,\dots,S_n;\tau) = \kappa(S_1,\dots,S_n).
\]
Using the result in the preceding section, we have

\medskip
\noindent
{\bf Fact.}
\[
\chi_{\mtop}(S_1,\dots,S_n) = \kappa(S_1,\dots,S_n) = \chi(S_1,\dots,S_n).
\]

\section{The Universal $n$-tuple of Self-adjoint Contractions}
\label{sec5}

In this section we prove
\[
\chi_{\mtop}(T_1,\dots,T_n) = \kappa(T_1,\dots,T_n)
\]
for the universal $n$-tuple of self-adjoint contractions $T_1,\dots,T_n$.

Let $A$ denote the unital $C^*$-algebra generated by $T_1,\dots,T_n$ or
equivalently the unital full free product of $n$ copies of $C([-1,1])$
with $T_j$ identified with the identical function in the $j$-th copy of
$C([-1,1])$.  $A$ has the property that given a $n$-tuple of self-adjoint
contractions $x_1,\dots,x_n$ in some unital $C^*$-algebra $B$, there is a
unique unital $*$-homomorphism of $A$ into $B$ mapping $T_j$ to $x_j$ $(1
\le j \le n)$.

Let $\tau(n)$ be the trace-state on $A$ which is the free product of $n$
copies of the trace-state $\tau(1)$ on $C([-1,1])$ which is given by the
equilibrium measure on $[-1,1]$, that is the arcsine distribution with
density $\pi^{-1}(1-x^2)^{-1/2}$ (\cite{4}).  We have
\[
\kappa(T_1,\dots,T_n) \le n\kappa(T_1) = n\chi(T_1;\tau(1)) =
\chi(T_1,\dots,T_n;\tau(n)).
\]

The $C^*$-algebra $A$ has sufficiently many finite-dimensional
representations.  This implies that given $\e > 0$ and noncommutative
polynomials $P_1,\dots,P_r$ we have
\[
\mvol(\G_{\mtop}(T_1,\dots,T_n;k_0,\e,P_1,\dots,P_r)) \ne 0
\]
for some $k_0 \in {\mathbb N}$.

It is easily seen that if $\d > 0$ is sufficiently small then for any
self-adjoint elements $x_1,\dots,x_n$ in a unital $C^*$-algebra so that
$\|x_j\| \le 1 + \d$ we will have $\|P_i(x_1,\dots,x_r)\| <
\|P_i(T_1,\dots,T_n)\| +
\e$
$(1
\le i \le r)$.  In particular using the usual microstates
\[
\G_{1+\d}(T_1,\dots,T_n;k,m,\e;\tau(n))
\]
for $T_1,\dots,T_n$ in $(A,\tau(n))$ we have
\begin{eqnarray*}
\G_{\mtop}(T_1,\dots,T_n;k_0,\e,P_1,\dots,P_r) &\oplus
&\G_{1+\d}(T_1,\dots,T_n;k,m,\e;\tau(n)) \\
&\subset
&\G_{\mtop}(T_1,\dots,T_n; k_0+k,\e,P_1,\dots,P_r)
\end{eqnarray*}
(the choice of parameters $m$, $\e$ and $\tau(n)$ in $\G_{1+\d}$ is
completely arbitrary).  This inclusion immediately yields
\[
\chi_{\mtop}(T_1,\dots,T_n) \ge \chi(T_1,\dots,T_n;\tau(n)).
\]

Putting things together {\em we conclude that we have proved}
\[
\chi_{\mtop}(T_1,\dots,T_n) = \kappa(T_1,\dots,T_n) =
\chi(T_1,\dots,T_n;\tau(n))
\]
{\em where $T_1,\dots,T_n$ is the universal $n$-tuple of self-adjoint
contradiction, and $\tau(n)$ the $n$-fold free product of the equilibrium
measure on $[-1,1]$}.

\section{Topological Free Entropy Dimension}
\label{sec6}

If $X$ is a subset of a metric space we denote by $N_{\e}(X)$ the minimum
number of elements in an $\e$-net of $X$.  The sets of norm-microstates
$\G_{\mtop}(\dots)$ will be viewed as subsets of $({\mathfrak M}_k)^n$
endowed with the uniform norm metric.  If $A,a_1,\dots,a_n,b_1,\dots,b_m$
are as in section \ref{sec2}, then we define
\[
D_{\e}(a_1,\dots,a_n: b_1,\dots,b_m)
\]
by taking
\[
\limsup_{k\rightarrow \i} k^{-2} \log N_{\e}(\G_{\mtop}(a_1,\dots,a_n:
b_1,\dots,b_m;k,\e,P_1,\dots,P_r))
\]
followed by an inf over $r \in {\mathbb N}$ and finite sets
$\{P_1,\dots,P_r\}$ of noncommutative polynomials.  Then the {\em
topological free entropy dimension of $a_1,\dots,a_n$ in the presence of
$b_1,\dots,b_m$ is defined by}
\[
\d_{\mtop}(a_1,\dots,a_n: b_1,\dots,b_m) = \limsup_{\e \rightarrow 0}
\frac {D_{\e}(a_1,\dots,a_n: b_1,\dots,b_m)}{|\log \e|}.
\]
Note that if we had used the normalized Hilbert--Schmidt norm
\[
|(A_1,\dots,A_n)|_2 = k^{-1/2}\left( \sum_{1 \le j \le n}
Tr_k(A_j^2)\right)^{1/2}
\]
instead of the uniform norm, we would have obtained the same quantity. 
Indeed,  a basic fact about $N_{\e}$ for unitary-invariant norms (\cite{6}
combined with Corollary~8 of \cite{5}) is that if $B(k,R)$ denotes the
ball of radius $R$ of ${\mathfrak M}_k^{sa}$ in uniform norm, then taking
$N_{\e}$ in uniform or normalized Hilbert--Schmidt norm we will have
\[
(C_1R/\e)^{k^2} \le N_{\e}(B(k,R)) \le (C_2R/\e)^{k^2}
\]
for some constants independent of $k$ and $\e \le R$.

Because of the upper bound for $N_{\e}(B(k,R))$ we always have
\[
\d_{\mtop}(a_1,\dots,a_n: b_1,\dots,b_m) \le n.
\]
To get $\d_{\mtop}$ for some $n$-tuples it is sufficient to remark the
following fact:  {\em if $\chi_{\mtop}(a_1,\dots,a_n: b_1,\dots,b_m) >
-\i$ then}
\[
\d_{\mtop}(a_1,\dots,a_n: b_1,\dots,b_m) = n.
\]
Indeed, we have
\[
N_{\e}(\G_{\mtop}(\dots))\mvol\ B^n(k,2\e) \ge \mvol(\G_{\mtop}(\dots))
\]
and by Corollary~4 of \cite{5} we have
\[
\liminf_{k \rightarrow \i} \left( \frac {n}{2} \log k + k^{-2} \log
\mvol(B^n(k,1))\right) > -\i
\]
which gives
\[
\liminf_{k \rightarrow \i} \left( \frac {n}{2} \log k + k^{-2} \log
\mvol(B^n(k,2\e))\right) \ge C + n \log \e.
\]
We infer that
\begin{eqnarray*}
\limsup_{k \rightarrow \i} k^{-2} \log N_{\e}(\G_{\mtop}(\dots)) &\ge &C +
n|\log \e| \\
& &+ \limsup_{k \rightarrow \i} \left(k^{-2} \log
\mvol(\G_{\mtop}(\cdots)) + \frac {n}{2} \log k \right)
\end{eqnarray*}
where $0 < \e < 1$.  This in turn immediately gives
$\d_{\mtop}(a_1,\dots,a_n: b_1,\dots,b_m) = n$ if
$\chi_{\mtop}(a_1,\dots,a_n: b_1,\dots,b_m) > -\i$.

It follows that {\em if $S_1,\dots,S_n$ is the semicircular $n$-tuple then}
\[
\d_{\mtop}(S_1,\dots,S_n) = n
\]
and similarly {\em if $T_1,\dots,T_n$ is the universal $n$-tuple of
self-adjoint contractions then}
\[
\d_{\mtop}(T_1,\dots,T_n) = n.
\]

A quantity with which $\d_{\mtop}$ should be compared is the {\em free
dimension capacity}
\[
\kappa\d(a_1,\dots,a_n) = \sup_{\tau \in TS(A)} \d_0(a_1,\dots,a_n;\tau)
\]
where $a_1,\dots,a_n$ is a self-adjoint generator of $A$ and
$\d_0(a_1,\dots,a_n;\tau)$ denotes the modified free-entropy dimension of
$a_1,\dots,a_n$ in $(A,\tau)$.  Clearly $\kappa\d(S_1,\dots,S_n) =
\kappa\d(T_1,\dots,T_n) = n$ by \cite{9}.

\section{Concluding Remarks}
\label{sec7}

\subsection{Semi-microstates}

A natural modification of the definition of $\G_{\mtop}$ is to consider
{\em norm-semi-microstates}, that is
\[
\G_{\mtop 1/2} (a_1,\dots,a_n: b_1,\dots,b_m;k,\e,P_1,\dots,P_r)
\]
with the inequalities
\[
|\|P_j(a_1,\dots,a_n,b_1,\dots,b_m)\| -
\|P_j(A_1,\dots,A_n,B_1,\dots,B_m)\|| \le \e
\]
replaced by
\[
\|P_j(A_1,\dots,A_n,B_1,\dots,B_m)\| \le
\|P_j(a_1,\dots,a_n,b_1,\dots,b_m)\| + \e.
\]

The corresponding quantity $\chi_{\mtop 1/2}$ clearly satisfies
$\chi_{\mtop 1/2} \ge \chi_{\mtop}$.  If for every $P_1,\dots,P_r,\e > 0$
there is some $k_0$ so that 
\[
\G_{\mtop}(a_1,\dots,a_n: b_1,\dots,b_m;k_0,\e,P_1,\dots,P_r) \ne \emptyset
\]
then actually
\[
\chi_{\mtop 1/2} = \chi_{\mtop}.
\]
Indeed, with $\oplus$ denoting componentwise orthogonal sums of operators,
we have
\begin{eqnarray*}
\G_{\mtop}(a_1,\dots,a_n: b_1,\dots,b_m;k_0,\e,P_1,\dots,P_r) &\oplus
&\G_{\mtop 1/2}(a_1,\dots,a_n: b_1,\dots,b_m;k,\e,P_1,\dots,P_r) \\
&\subset &\G_{\mtop}(a_1,\dots,a_n: b_1,\dots,b_m;k+k_0,\e,P_1,\dots,P_r).
\end{eqnarray*}

\subsection{Orthogonal Sums}

Let $\a' = (a'_1,\dots,a'_n)$, $\a'' = (b'_1,\dots,b'_m)$, $\a' \oplus
\a'' = (a'_1 \oplus a''_1,\dots,a'_n \oplus a''_n)$, $\b' =
(b'_1,\dots,b'_m)$ etc.  It is easily seen that
\begin{eqnarray*}
\G_{\mtop}(\a': \b';k',\e,P_1,\dots,P_r) &\oplus &\G_{\mtop}(\a'':
\b'';k'',\e,P_1,\dots,P_r) \\
&\subset &\G_{\mtop}(\a' \oplus \a'': \b' \oplus \b'';k' +
k'',\e,P_1,\dots,P_r)
\end{eqnarray*}
so that
\[
\chi_{\mtop}(\a' \oplus \a'': \b' \oplus \b'') \ge \max(\chi_{\mtop}(\a':
\b'),\chi_{\mtop}(\a'': \b'')).
\]
On the other hand it is easily seen that
\[
\kappa(\a' \oplus \a'': \b' \oplus \b'') = \max(\kappa(\a':
\b'),\kappa(\a'': \b'')),
\]
since the extremal trace-states of $C^*(\a' \oplus \a'', \b' \oplus \b'')$
are either trace-states of $C^*(\a',\b')$ or of $C^*(\a'',\b'')$.

Therefore a test for the general question about the equality of
$\chi_{\mtop}$ and $\kappa$ is whether we actually have
\[
\chi_{\mtop}(\a' \oplus \a'': \b' \oplus \b'') = \max(\chi_{\mtop}(\a':
\b'),\chi_{\mtop}(\a'': \b'')).
\]

\subsection{Free Equilibrium Trace-States}

In classical potential theory, the probability measures minimizing the
logarithmic energy are the equilibrium measures.  Generalizing this from
the one-variable case, we arrive at the question of finding those $\tau
\in TS(A)$ for which $\kappa(a_1,\dots,a_n) = \chi(a_1,\dots,a_n;\tau)$
(the generator $a_1,\dots,a_n$ of $A$ is given).

Both trying to find general facts about such $\tau$, as well as finding
such $\tau$ for specific $a_1,\dots,a_n$ seem to be questions which deserve
some attention.

\subsection{Topological Free Entropy Dimension of Generators}

By analogy with the von~Neumann algebra case, it seems natural to ask
whether the topological free entropy dimension of a generator is an
invariant of the $C^*$-algebra; i.e., if $a'_1,\dots,a'_n$ and
$a''_1,\dots,a''_m$ are self-adjoint generators of the $C^*$-algebra $A$,
does it follow that $\d_{\mtop}(a'_1,\dots,a'_n)$ and
$\d_{\mtop}(a''_1,\dots,a''_m)$ are equal?

\end{document}